\documentclass[10pt]{article}

\newtheorem{theorem}{Theorem}[section]

\newtheorem{definition}{Definition}

\usepackage{epsfig} 

\title{Pontryagin's principle of stabilization}

\author{Sergey Nikitin \\
\thanks{Department of Mathematics and Statistics, Arizona State University, Tempe, AZ 85287-1804  {\tt \small nikitin@asu.edu}}
}
 
\begin{document}

\maketitle

\begin{abstract}

The paper presents necessary and sufficient conditions for a nonlinear system to be stabilized by a feedback. The conditions are based on the ideas related to the well-known Pontryagin's maximum principle. That allows us to formulate the results in terms that are valid for continuous, discontinuous, stationary and time-dependent feedbacks.
\end{abstract}

\section{INTRODUCTION}

Stabilization is one of the central topics of control theory. It was shown in \cite{Aeyels}, \cite{Coron}, \cite{Ser_book}, \cite{Ser_nec_stab} that the class of continuous stationary feedbacks is too restrictive for the purposes of stabilization of nonlinear systems. In other words, in order to design a continuous feedback stabilizer one needs to use functions depending on time and state \cite{pomet}, \cite{samson}. To stay with the class of stationary feedbacks one has to deal with piecewise continuous functions \cite{Clark}, \cite{Nikitin99}.  The synthesis procedures for both continuous and piecewise continuous stabilizers are developed only for some special types of systems. For example, in  \cite{Nikitin99}, \cite{Ser_book} it is shown how to construct stationary piecewise continuous stabilizers for generic two-dimensional affine nonlinear systems. On the other hand, the papers \cite{pomet}, \cite{sordalen}, \cite{samson} show how to design feedbacks for certain types of nonholonomic systems.  Both approaches (nonstationary continuous, and stationary discontinuous) are quite complicated as far as the feedback synthesis is concerned. 

\vspace{0.1cm}

 This paper presents necessary and sufficient conditions for a nonlinear system to be stabilized by a  feedback. The main results are presented in terms that are valid for continuous, discontinuous, stationary and time dependent feedbacks. Our approach is based on the properties of the Pontryagin type Hamiltonian systems. The proposed sufficient conditions lead us to an effective feedback synthesis procedure that allows to construct piece-wise continuous stabilizing feedback laws for a general nonlinear system.

\vspace{0.1cm}

\section{PRELIMINARIES}
Consider a system
\begin{equation}
\label{system}
\dot x =  f(t,x, u)
\end{equation}
where $u$ is the control input; $x$ denotes the state of the system and $x\in {\rm R}^n,$ $ {\rm R}^n$ -- $n$-dimensional linear real space. $f(t,x,u) $ is a vector field:  
$$
\forall (t,x, u) \in {\rm R} \times {\rm R}^{n} \times {\rm R}^m\;\;\;f(t,x,u)\in {\rm R}^n
$$
and $f(t,x,u)$ is continuously differentiable with respect to $(t,x,u).$ We write $f\in {\rm C^1}.$
Throughout the paper we assume that  ${\rm R}^n$ is equipped with the scalar product and $\| x \|$ denotes the magnitude of $x,$ i.e
$\| x \| =\sqrt {\langle x , x \rangle },$ where $\langle x , x \rangle$ is the scalar product of $x$ with itself.

Consider the initial value problem
\begin{eqnarray}
\label{initialValueProblem}
\dot x (t) &=&  f(t,x, u(t,x))  ,\nonumber\\
&&\\
x(t_0)&=&x_0, \nonumber
\end{eqnarray}
where $f \in {\rm C}^1 $ and $u(t,x)$ is the feedback. 
Its solution is defined in Filippov's sense \cite{Filippov}.

\vspace{0.1cm}

Let $x=0$ be the Filippov solution for (\ref{system}) where $t\in {\rm R}$ and $u=u(t,x)$ is a feedback law. We call such solution an equilibrium of (\ref{system}).
The goal of this paper is to find necessary and sufficient conditions for the feedback
\begin{equation}
\label{feedback}
u= u(t,x) 
\end{equation}
to stabilize the system (\ref{system}) at the origin over a compact set  $K \subset  {\rm R}^n.$  Throughout the paper we assume that $u(t,x)$ takes its values from a compact subset 
\begin{equation}
\label{controlConstr}
\Omega \subseteq {\rm R}^m .
\end{equation}
The stabilization is defined as follows.

\begin{definition}
\label{stabilizationDef} 
The system (\ref{system}) is said to be stabilizable over a compact set $K \subset  {\rm R}^n$  by the feedback (\ref{feedback}) if the solutions $x(t,t_0,x_0)$ of the closed loop system (\ref{initialValueProblem}) satisfy the condition
$$
\lim_{t \to \infty } x(t,t_0,x_0) = 0  \;\;\forall\; x_0 \in K , \;\;t_0\in {\rm R}  .
$$
and the equilibrium is stable.
\end{definition}

\section{NECESSARY AND SUFFICIENT CONDITIONS}

Our main goal is to find necessary and sufficient conditions for the feedback (\ref{feedback}) to stabilize the system (\ref{system}). Consider the optimal control problem (with free end)

\begin{eqnarray}
\label{optimal1}
\int_{t_0}^T \frac{\partial}{\partial t} V(t,x(t,t_0,x_0)) +&& \nonumber \\
&&\\
 \langle \frac{\partial }{\partial x} V(t,x(t,t_0,x_0)) , f(t,x(t,t_0,x_0),u(t)) \rangle dt  &\to& \inf_{u(t) },\nonumber
\end{eqnarray} 
where $V(t,x)$ is a Lyapunov type function (see \cite{Rouche}, \cite{Zubov})  and  $x(t,t_0,x_0)$ is the solution for 
\begin{eqnarray}
\label{optimal2}
\dot x(t) &=& f(t,x,u), \nonumber\\
&& \\
x(t_0) &=& x_0 , \nonumber
\end{eqnarray}
where the control input $u$ takes its values from the set $\Omega $ defined in (\ref{controlConstr}).

The existence of the optimal solution for (\ref{optimal1}), (\ref{optimal2}) was studied in great details \cite{Berkovitz}. For simplicity we assume that
$$
\frac{\partial }{\partial x} f(t,x,u) \;\;\mbox{ and } \;\;\frac{\partial }{\partial u} f(t,0,u)
$$
are uniformly bounded for $(t,x,u) \in {\rm R} \times {\rm R}^n \times \Omega : \;\;\exists \; Q>0$ such that $\forall \;\;t\in {\rm R}, \;x\in {\rm R}^n, \; u \in \Omega\;$
\begin{eqnarray}
\label{bound}
\|\frac{\partial }{\partial x} f(t,x,u) x \|&\le& Q\cdot \|x\| \nonumber\\
 \mbox{ and } &&\\
 \|\frac{\partial }{\partial u} f(t,0,u) u \| &\le& Q\cdot \|u\|. \nonumber
\end{eqnarray}

Then the existence of the optimal trajectory for any $x_0$  follows from the classical Arzela's theorem. The conditions (\ref{bound}) can be relaxed (at the expense of the simplicity of formulations) with the help of the results from \cite{Berkovitz}

\vspace{0.1cm}

The necessary conditions for optimality are provided by Pontryagin result \cite{Pontryagin}. Let
\begin{eqnarray}
H(\psi ,t, x , u ) = \frac{\partial}{\partial t} V(t,x) &+& \nonumber \\
 \langle \frac{\partial }{\partial x} V(t,x) , f(t,x,u) \rangle &+& \langle \psi , f(x,u) \rangle . \nonumber
\end{eqnarray}
Then the optimal solution satisfies
\begin{eqnarray}
\dot x &=& \frac{\partial }{\partial \psi } H(\psi ,t, x , u ), \nonumber \\
\dot \psi &=& - \frac{\partial }{\partial x } H(\psi ,t, x , u ), \nonumber
\end{eqnarray}
where $\psi (T)=0$ and $u=u(t,x,\psi)$ is the solution for
$$
    H(\psi ,t, x , u ) \to \inf_{u \in \Omega } .
$$

 Introduce notations
$$
\nu =   \frac{\partial }{\partial x} V(t,x) + \psi 
$$
and
$$
S(\nu ,t, x ,u) = \langle \nu , f(t,x,u) \rangle .
$$
Then the Pontryagin necessary conditions of optimality take the following form
\begin{eqnarray}
\label{Pontryagin1}
\dot x &=& \frac{\partial }{\partial \nu } S(\nu ,t,x , u ), \nonumber \\
&&\\
\dot \nu &=& - \frac{\partial }{\partial x } S(\nu ,t,x , u ), \nonumber
\end{eqnarray}
where $\nu (T) =  \frac{\partial }{\partial x} V(T,x(T,x_0))$ and $u=u(t,x,\nu)$ is the solution for
\begin{equation}
\label{Pontryagin2}
S(\nu ,t, x ,u) \to \inf_{u \in \Omega }.
\end{equation}

Now we can formulate our necessary conditions of stabilization.
\begin{theorem}
\label{necess1} 
If the system (\ref{system}) with $f\in C^1$ and (\ref{bound})  is stabilizable over $K,$ then 
$$
\forall \;\; \delta > 0,\;\;t_0\in {\rm R},\;\; x_0 \in \; K \;\;
$$
$$
\exists \;\; \nu_0 \in {\rm R}^n,\;\; \mbox{ and } T >0
$$
such that for the solution $(x(t), \;\nu(t) )$ of (\ref{Pontryagin1}), (\ref{Pontryagin2}) with the initial conditions 
$$
x(t_0) = x_0,\;\;\nu(t_0) = \nu_0, 
$$
we have
$$
\| x(T) \| \le \delta  .
$$
and
$$
 \nu(T) =  x(T).
$$
\end{theorem}

\vspace{0.1cm}

{\bf Proof.}
Let the system (\ref{system}) be stabilizable by a feedback $u(t,x)$ over a compact set $K.$  Then 
\begin{equation}
\label{stb}
\forall t_0\in{\rm R},  \;\;x_0\;\in\;K\;\;\lim_{t \to \infty } x_u(t,t_0,x_0) = 0,
\end{equation}
where $x_u(t,t_0,x_0)$ denotes the solution of the initial value problem
\begin{eqnarray}
\dot x(t) &=& f(t,x,u(t,x)), \nonumber\\
x(t_0) &=& x_0 . \nonumber
\end{eqnarray}
We need to prove that for any $\delta > 0,\;\;t_0\in {\rm R} $ and any $x_0\;\in\;K\subseteq {\rm R}^n$ there exist $\nu_0 \;\in\;{\rm R}^n$ and $T>0$ such that 
$$
 \|x(T)\| < \delta \;\;\mbox{ and } \nu (T) = x(T)\;\;
$$
where $x(t)$ is the $x$-part of the solution  $(x(t), \;\nu(t) )$ of (\ref{Pontryagin1}), (\ref{Pontryagin2}). 

\vspace{0.1cm}

 Consider the extremal problem (\ref{optimal1}), (\ref{optimal2}), where $V(t,x) = \frac12 \|x\|^2.$ It follows from (\ref{bound}) and the classical Arzela's theorem that the optimal solution exists. Hence, it satisfies Pontryagin necessary conditions, 

$$
\|x(T)\| \le \|x_u(T,t_0,x_0)\|, \;\;\;\;\nu (T) =  x(T)
$$
and (\ref{stb}) yields the existence of $T>0$ such that
$$
\| x(T) \| \le \delta \;\;\mbox{ and } \nu (T) =x(T) .
$$
{\bf Q.E.D.}

\vspace{0.1cm}

Let $u=u(t,x,\nu)$ denote the control such that
$$
\langle \nu , f(t,x,u(t,x,\nu)) \rangle = \inf_{u\in \Omega } \langle \nu , f(t,x,u(t,x,\nu)) \rangle .
$$
Then it follows from 
$$
\langle \nu , f(t,x,u(t,x,\nu)) \rangle \le \langle \nu , f(t,x,u) \rangle  \;\;\forall\;u\in \Omega
$$
that for any $\lambda > 0$
$$
\langle \lambda \nu , f(t,x,u(t,x,\nu)) \rangle \le \langle \lambda \nu , f(t,x,u) \rangle  \;\;\forall\;u\in \Omega.
$$
Hence, for $\lambda >0$ $u(t,x,\lambda \cdot \nu) = u(t,x,\nu)$ and the Hamiltonian $S(t,x,\nu)=S(\nu, t,x,u(t,x,\nu))$ has the following important property
\begin{equation}
\label{Guigens}
S(t,x, \lambda \nu ) = \lambda S(t,x, \nu )\;\;\mbox{ for }\;\; \lambda >0.
\end{equation}
Thus, 
$$
      \nu \cdot dx = 0
$$
is preserved along the solutions of the system (\ref{Pontryagin1}), (\ref{Pontryagin2}). Consider the $(n-1)$-dimensional manifold $M_\delta$ defined as
$$
\nu = grad V(x)\;\;\mbox{ and } V(x) = \delta, 
$$
where $\delta > 0$ and $V(x)$ is a Lyapunov function. Clearly, restriction of the differential form $\nu \cdot dx$ to this manifold is equal to zero. The Lagrangian manifold
$L_\delta$ defined by the trajectories of  (\ref{Pontryagin1}), (\ref{Pontryagin2}) emitted from  $M_\delta$ plays an important role in analyzing stabilization and in designing feedback stabilizers. In order to illustrate the significance of $L_\delta $ for the stabilization problems we turn our attention to stationary control systems having the form
$$
\dot x = f(x,u).
$$
The Lagrangian manifold $L_\delta$ we define as the set of trajectories for 
\begin{eqnarray}
\label{Pontryagin3}
&&S(\nu ,x ,u) \to \inf_{u \in \Omega } \;\nonumber \\
\dot x &=& -\frac{\partial }{\partial \nu } S(\nu ,x , u ), \nonumber \\
&&\\
\dot \nu &=&  \frac{\partial }{\partial x } S(\nu ,x , u ), \nonumber
\end{eqnarray}
emitted from the manifold $M_\delta ,$ where
$$
 S(\nu ,x , u ) = \langle \nu , f(x,u) \rangle.
$$

 Taking a fixed point $x_0\in M_\delta$ we obtain the function
\begin{equation}
\label{generating}
W(P) =  \int_{x_0}^{P} \nu dx + V(x_0)
\end{equation}
$W(P)$ is well defined on $L_\delta .$

 The property (\ref{Guigens}) implies the Guigens principle for the propagation of "light" described by (\ref{Pontryagin3}). In other words, the property (\ref{Guigens}) suggests (see \cite{Maslov}) that we can treat stabilization similar to a problem from optics.  In the terminology of optics Theorem \ref{necess1} states that for the system to be stabilizable at the origin  over a compact set $K$ it is necessary that the sources of "light" located on the surface of any tiny ball centered at the origin  illuminate all points from $K.$  The tiny ball can be replaced by the set where $V(x)\le \varepsilon.$ Moreover, the wave fronts are defined by the projections to $x$-space the level sets of the function $W(P).$ Under certain conditions (see \cite{Nikitin99})  Theorem \ref{necess1} is not only necessary but also sufficient for stabilization. Moreover, the Pontryagin's principle hints the form of the feedback stabilizer outside the set where $V(x)\le \varepsilon.$ In other words, the following statement takes place. 

\begin{theorem}
\label{design}
Let the system 
$$
\dot x = f(x,u)
$$
be locally stabilizable at the equilibrium $x=0.$ Moreover, there exists a stabilizing feedback law $u=w(x)$ and a Lyapunov function $V(x)$ such that for sufficiently small $\varepsilon > 0$ we have
$$
\langle grad  V(x) , f(x,w(x))\rangle \le 0\;\;\mbox{ for } \;\;\; 0< V(x)\le \varepsilon .
$$ 
If the Lagrange manifold $L_\varepsilon$ is such that the restriction of $u(\nu ,x )$ (from (\ref{Pontryagin3})) to $L_\varepsilon$ is a function of $x$ then the system is stabilizable over any set $K$ which is covered by the natural projection of $L_\varepsilon $ onto the $x$-space.
\end{theorem}

\vspace{0.1cm}

{\bf Proof.} We follow the ideas presented in \cite{Nikitin99} and choose the feedback $u=u(x)$ equal to $w(x)$ for all $x$ such that $V(x)\le \varepsilon.$ On the other hand, if $V(x) > \varepsilon$ and $x$ belongs to the projection of $L_\varepsilon$ onto the $x$-space then we define the feedback as
$$
  u=u(\nu,x),
$$
where the function $u(\nu,x)$ is defined by (\ref{Pontryagin3}) along the lines of well-known Pontryagin's procedure and $ u=u(\nu,x)$ is a function of $x$ by the assumption. 

\vspace{0.1cm}

{\bf Q.E.D.}

\vspace{0.1cm}

Theorem \ref{design} suggests an efficient feedback design procedure. However, it tacitly imposes some restrictions in the form of requirements fulfilled by certain properties of the Lagrangian manifold $L_\varepsilon.$ Although the stabilization problem can be reformulated in optical terms the machinery of Lagrangian manifolds and symplectic geometry is not readily available for the control problems. The main obstacle here lies in the lack of smoothness of the corresponding Hamiltonian. However, since the property (\ref{Guigens}) implies the Guigens principle for the propagation of "light" described by (\ref{Pontryagin3}) the Lagrangian manifold $L_\varepsilon$ can be constructed numerically. Moreover, if the Lagrangian manifold $L_\varepsilon$ does not admit the parametrization by $x$ then the corresponding feedback can be designed (for some applications) with the ideology provided by the theory of viscosity solutions (see, e.g., \cite{Day}). 

Theorem \ref{design} leads us to an effective design synthesis of feedback laws suitable for various applications of control theory. The format and the size of this brief note does not allow us to consider in details various applied versions of Theorem \ref{design}. We only illustrate the power of Theorem \ref{design} by the results on stabilization of an affine nonlinear system with one control input,
$$
\Sigma(f,b):\;\;\;\dot x = f(x) + b(x) u,
$$
where $x\in {\rm R}^n$ and $f(x), b(x)$ are $C^\infty$ vector fields on ${\rm R}^n.$ We assume that $f(0)=0$ and the control $u$ takes its values from $[-1,\; 1].$ Then (\ref{Pontryagin3}) yields 
\begin{eqnarray}
\label{Pontryagin4}
\dot x &=& -f(x)+sgn(\langle \nu, b(x) \rangle ) \cdot b(x), \nonumber \\
&&\\
\dot \nu &=&  (\frac{\partial }{\partial x }f(x))^T\nu - sgn(\langle \nu, b(x) \rangle ) \cdot  (\frac{\partial }{\partial x }b(x))^T \nu , \nonumber
\end{eqnarray}
where $sgn(\langle \nu, b(x) \rangle )$ denotes the sign of $\langle \nu, b(x) \rangle.$
\begin{theorem}
\label{affineSys}
Let the system $\Sigma(f,b)$ be locally stabilizable at the equilibrium $x=0.$ Moreover, there exist a stabilizing feedback law $u=w(x)$ and a Lyapunov function $V(x)$ such that for sufficiently small $\varepsilon > 0$ we have
$$
\langle grad  V(x) , f(x) + b(x)\cdot w(x)\rangle \le 0\;\; \mbox{and}\;\;\vert w(x)\vert \le 1 
$$
when $ 0< V(x)\le \varepsilon.$
Assume 
\begin{equation}
\label{transversal}
\langle grad  V(x) , f(x)+b(x)\cdot w(x)\rangle  < 0 \;\;\;\mbox{ for } \;\;\; V(x) =  \varepsilon 
\end{equation}
and
\begin{equation}
\label{transvers}
<\nu, ad_f b(x)>\not=0.
\end{equation}
 for any point $(x,\nu) $ from the Lagrange manifold $L_\varepsilon$ such that
\begin{equation}
\label{switching}
\langle \nu, b(x) \rangle =0.
\end{equation}
Then the system is stabilizable at the origin over any set $K$ which is covered by the natural projection of $L_\varepsilon $ onto the $x$-space. 
\end{theorem}

{\bf Proof.} Consider the coordinate charts $(\psi, \tau )$  on the Lagrangian manifold $L_\varepsilon,$  where $\psi$ denotes local coordinates on the surface $$
V(x)=\varepsilon
$$
and $\tau$ is the parameter along the solutions of (\ref{Pontryagin4}) known as bicharacteristics (see. e.g., \cite{Maslov} for details and further references).
Due to the special form of (\ref{Pontryagin4}) and the condition (\ref{transversal}) we obtain with the help of Liouville theorem that the Jacobian 
$$
det(\frac{\partial x}{\partial \psi }, \frac{\partial x} {\partial \tau })
$$
is not equal to zero along the bicharacteristics. Therefore, the Lagrange manifold $L_\varepsilon $  locally admits a generating function $W(x)$ (defined as (\ref{generating}), see also \cite{Maslov_Fedoriuk}). Hence, the control 
$$
v(x)=-sgn(\langle \nu, b(x) \rangle )= -sgn(\langle \frac{\partial W}{\partial x}, b(x) \rangle )
$$
depends only on $x$ as long as the solutions of (\ref{Pontryagin4}) are transverse to the switching surface defined by (\ref{switching}). The latter is assured by (\ref{transvers}). 

\vspace{0.1cm}

{\bf Q.E.D.}

\vspace{0.1cm}

\begin{figure}
      \centering
\epsfxsize=9cm
\epsfysize=7cm
\hbox{\centerline{\epsfbox{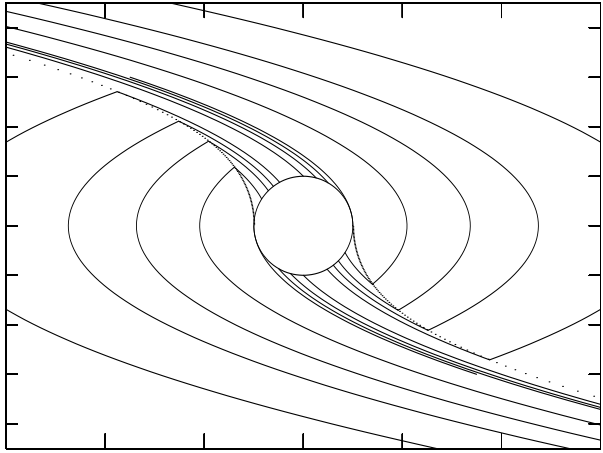}}}
      \caption{ The "light" emitted from the points of the circle illuminates all points outside the disk.}
      \label{lightFromCircle}
   \end{figure}

To illustrate an application of the Pontryagin's stabilization principle presented in this paper consider the  system 
\begin{eqnarray}
\label{chain}
\dot x_1 &=& x_2, \nonumber \\ 
\dot x_2 &=& u. \nonumber 
\end{eqnarray}
We show that for any number $C>0$ this system can be globally stabilized at the origin by a piece-wise smooth feedback
$$
u=u(x) \mbox{ such that } \vert u(x)\vert \le C \;\;\forall \; x\in {\rm R}^2.
$$ 
Consider the Pontryagin Hamiltonian
$$
S(x,\nu) =  \nu_1 \cdot x_2  - \vert \nu_2 \vert.
$$
The system  corresponding to (\ref{Pontryagin3}) takes the form
\begin{eqnarray*}
\dot x_1 &=& - x_2, \nonumber \\ 
\dot x_2 &=& sgn(\nu_2),  \nonumber
\end{eqnarray*}
and
\begin{eqnarray*}
\dot \nu_1 &=& 0, \nonumber \\ 
\dot \nu_2 &=& \nu_1 \nonumber
\end{eqnarray*}
where $sgn(z)$ denotes the sign of $z.$ The "light" emitted from the points of the circle of radius $\delta$ centered at the origin illuminates all points outside the disk (Fig.\ref{lightFromCircle}).  

\begin{figure}
      \centering
\epsfxsize=9cm
\epsfysize=7cm
\hbox{\centerline{\epsfbox{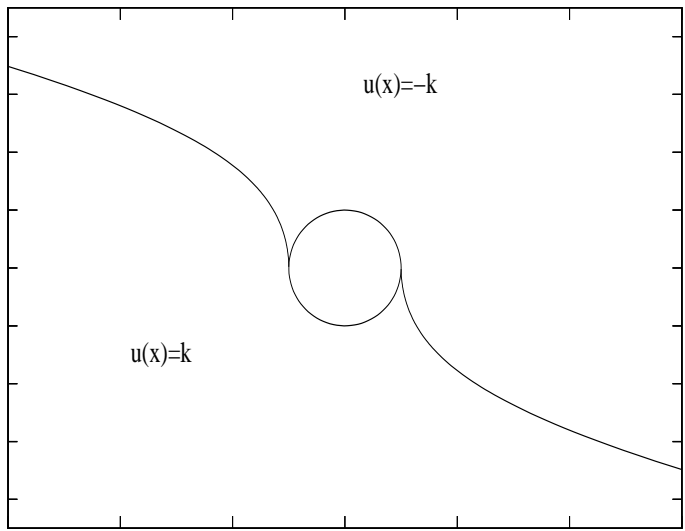}}}
      \caption{The feedback design outside the disk.}
      \label{switch}
   \end{figure}

 Moreover, the feedback
$$
w(x) = -x_1 -x_2
$$ 
stabilizes the system over the disk of radius $\delta$ and the disk is invariant with respect to the closed loop system. Outside the ball we define a piece-wise constant control as shown in Fig.\ref{switch}. The constant $\delta$ as well as control values $-k,\;k$ can be chosen so that 
$$
\vert u(x) \vert \le C \;\;\forall \;x \in {\rm R}^2.
$$
Notice, that for $\delta=1$ the switching curve (the curve where the control switches between the values $-k$ and $k)$ is defined by
\begin{eqnarray*}
x_1 &=& -\frac{\sin(\tau )\vert \sin(\tau)\vert }{2\cdot \cos^2(\tau )} + \frac{\sin^2 (\tau )}{\cos(\tau )} + \cos(\tau ) \\
x_2 &=& - \frac{\vert \sin(\tau )\vert }{\cos (\tau )} + \sin (\tau ),
\end{eqnarray*}
where $\frac{\pi}{2} < \tau < \pi$ and $\frac{3 \pi }{2} < \tau < 2\cdot \pi .$

Consider the affine nonlinear system having the form
$$
\Sigma(f,b):\;\;\;\dot x = f(x) + b(x) u,
$$
where $f(x),\;\;b(x): {\rm R}^2 \to {\rm R}^2 $ are twice continuously differentiable mappings. Following the scenario from Chapter 2 of \cite{Ser_book} we introduce the equilibrium set
$$
\varphi^{-1} (0) = \{ x \in {\rm R}^2;\;\;det (f(x), b(x)) = 0  \}.
$$
Making use of the theory developed in \cite{Ser_book} (see Chapter 2) we come to the following result.

\begin{theorem}
\label{2D}
If
\begin{equation}
\label{fullrank}
rank \{ b(x) , ad_f b(x) \} = 2 \;\;\;\forall\;x\in {\rm R}^2 
\end{equation} 
then the bang-bang controllability of $\Sigma(f,b)$ implies its global (semi-global) stabilization at any equilibrium from $\varphi^{-1}(0)$  by a piece-wise smooth feedback $u(x)$ such that
$$
 \vert u(x) \vert \le 1\;\;\forall\;x\in {\rm R}^2.
$$
\end{theorem}

{\bf Proof.} It follows from (\ref{fullrank}) (see Chapter 2 from \cite{Ser_book}) that for any point from the equilibrium set one can construct a locally stabilizing feedback and the corresponding Lyapunov function. Moreover, 
$$
 \langle \nu ,ad_f b \rangle \not=0 
$$
when $ \langle \nu , b \rangle =0 .$
 Hence, by Theorem \ref{affineSys} there exists the stabilizing feedback for system $\Sigma(f,b).$  

 \vspace{0.1cm}

{\bf Q.E.D.}

\vspace{0.1cm}

Theorem \ref{2D} leads us to an efficient controller design procedure for manipulators. In order to illustrate its application consider a manipulator defined by 
\begin{eqnarray}
\label{manipulator}
\dot x_1 &=& x_2 \nonumber\\
&&\\
\dot x_2 &=& f(x_1,x_2) + u, \nonumber
\end{eqnarray}
where $f(x_1,x_2)$ is a smooth function such that 
$$
f(0,0)=0.
$$
Suppose we can measure only $x_1(t)$ and $x_2(t)$ is not available for us. Then we design the stabilizing feedback law in the following two steps.
\begin{itemize}
\item[i.] Design the feedback $u(x_1,x_2)$ with the help of Pontryagin's principle.
\item[ii.] Construct an estimator for $x_2$ and replace $u(x_1,x_2)$ with $u(x_1,z_2),$ where $z_2$ is the corresponding estimation for $x_2(t).$
\end{itemize}
The first step has been already described in details in this paper. Now we outline the construction of the estimator mentioned in the second step and then show that the resulted feedback law delivers stabilization for the manipulator (\ref{manipulator}) under some generic conditions.

The estimator is given by
\begin{eqnarray}
\label{estimator}
\dot z_1 &=& z_2 - \beta_1(z_1 -x_1 ) \nonumber\\
&&\\
\dot z_2 &=& f(x_1,z_2) + u(x_1,z_2) - \beta_2(z_1 -x_1 ), \nonumber
\end{eqnarray} 
where $\beta_1,\;\;\beta_2$ are positive real numbers. Let $W(x_1,x_2)$ denote the generating function (\ref{generating}) corresponding to the Pontryagin's type control $u(x_1,x_2).$ Then the following statement takes place.

\begin{theorem}
\label{manipulatorFeedback}
If the system (\ref{manipulator}) is bang-bang controllable and 
$$
f(x_1,x_2),\;\;\frac{\partial}{\partial x_2} W(x_1,x_2)
$$
 are globally Lipschitz with respect to the second argument then there exist real positive numbers $\beta_1,\;\;\beta_2$ such that the Pontryagin's type feedback law $u(x_1,z_2),$ where $z_2$ is calculated by (\ref{estimator}) stabilizes the system (\ref{manipulator}) at the origin. 

\end{theorem}
{\bf Proof.} By Theorem \ref{2D} there exists a Pontryagin's type feedback $u(x_1,x_2)$ that stabilizes the system at the origin. Since $x_2(t)$ is not available for us we need to use estimator (\ref{estimator}). Consider a Lyapunov function with $e_1=z_1-x_1$ and $e_2=z_2-x_2$:
$$
V=2\frac{\beta_2}{\beta_1} e_1^2 - 2 e_1\cdot e_2 + (\frac{2}{\beta_1} + \frac{\beta_1}{\beta_2})e_2^2,
$$
where $\beta_1,\;\beta_2 $ are positive numbers. Differentiating this function with respect to the closed loop system $(u= u(x_1,z_2))$  (\ref{manipulator}), (\ref{estimator}) yields
$$
\frac{d}{dt} V = - 2\beta_2 \cdot e_1^2 - 2 e_2^2 +
$$
$$
2\cdot(-e_1 + (\frac{2}{\beta_1} + \frac{\beta_1}{\beta_2})e_2)(f(x_1,z_2)-f(x_1,x_2))
$$
 \vspace{0.1cm}
Taking into account that 
$$
\vert f(x_1,z_2)-f(x_1,x_2) \vert \le L \cdot \vert z_2 - x_2 \vert \;\;\forall\;x_1,\;x_2,\;z_2
$$
and
$$
2 \vert e_1 \cdot e_2 \vert \le (\frac{e_1}{\delta})^2 + (\delta \cdot e_2)^2
$$
we obtain that 
$$
\frac{d}{dt} V  \le ( - 2\beta_2  + \frac{1}{\delta^2}\cdot L) \cdot e_1^2 +
$$
$$
(-2 + \delta^2 \cdot L + (\frac{2}{\beta_1} + \frac{\beta_1}{\beta_2})\cdot L) \cdot e_2^2
$$
It is evident that for any $L>0$ there exist $\delta>0$ and positive $\beta_1,\;\;\beta_2$ such that 
$$
- 2\beta_2  + \frac{1}{\delta^2}\cdot L <0 \;\;
$$
and
$$
-2 + \delta^2 \cdot L + (\frac{2}{\beta_1} + \frac{\beta_1}{\beta_2})\cdot L <0.
$$
Hence, (see, e.g., \cite{Queiroz}) 
$$
e_1(t) \to 0\;\;\mbox{ and } \;\;e_2(t) \to 0\;\;\;\mbox{ as } t\to \infty.
$$
Now let us differentiate $W(x_1,x_2)$ with respect to the closed loop system (\ref{manipulator}).
$$
\frac{d}{dt} W = \nu_1 \cdot x_2 + \nu_2(f(x_1,x_2) + u(x_1,x_2)) +
$$
$$
 \nu_2(u(x_1,z_2) - u(x_1,x_2)),
$$
where 
$$
(\nu_1,\;\nu_2)= (\frac{\partial}{\partial x_1} W,\; \frac{\partial}{\partial x_2} W).
$$
By construction of the Pontryagin's type feedback there exist a positive number $\gamma > 0$ such that 
$$
\nu_1 \cdot x_2 + \nu_2(f(x_1,x_2) + u(x_1,x_2)) < -\gamma
$$
for all $(\nu,\;x)$ from the corresponding Lagrangian manifold. $\frac{\partial}{\partial x_2} W(x_1,x_2)$ is globally Lipschitz with respect to the second argument: $\forall\;\;x_1,\;x_2,\;z_2$ 
$$
 \vert \frac{\partial}{\partial x_2} W(x_1,x_2)  -   \frac{\partial}{\partial x_2} W(x_1,z_2) \vert \le M \cdot \vert x_2 - z_2 \vert.
$$
Taking into account
$$
u(x_1,x_2) = -sgn(\frac{\partial}{\partial x_2} W(x_1,x_2))
$$
we have
$$
\frac{\partial}{\partial x_2} W(x_1,x_2) \cdot (-sgn(\frac{\partial}{\partial x_2} W(x_1,z_2))+
$$
$$
 sgn(\frac{\partial}{\partial x_2} W(x_1,x_2))) =
$$
$$
 (\frac{\partial}{\partial x_2} W(x_1,x_2) - \frac{\partial}{\partial x_2} W(x_1,z_2))\cdot (-sgn(\frac{\partial}{\partial x_2} W(x_1,z_2))) +
$$
$$
 \vert \frac{\partial}{\partial x_2} W(x_1,x_2) \vert -  \vert \frac{\partial}{\partial x_2} W(x_1,z_2) \vert.
$$
Hence,
$$
\vert \nu_2(u(x_1,z_2) - u(x_1,x_2)) \vert \le 2 M \vert e_2(t) \vert.
$$
Since $e_2(t) \to 0$ as $t\to \infty$ one can find $T>0$ such that 
$$
\vert \nu_2(u(x_1,z_2) - u(x_1,x_2)) \vert \le 2 M \vert e_2(t) \vert < \frac{\gamma}{2} \;\;\forall\;t> T.
$$
Thus, 
$$
\frac{d}{dt}W <-\frac{\gamma}{2}\;\;\forall\;t> T
$$
and at some moment of time the solution will arrive at the neighborhood of the origin where one can switch to the feedback control constructed with the help of the classical linear theory.

{\bf Q.E.D.}



\end{document}